\documentclass[10pt,reqno]{article}
\usepackage{latexsym, amsmath, amssymb}
\usepackage{graphicx,epsfig}

\usepackage{amsfonts}
\usepackage{amscd}
\usepackage{amsbsy}

\newtheorem{thm}{Theorem}[section]

\newtheorem{lem}[thm]{Lemma}

\numberwithin{equation}{section}
\newcounter{saveeqn}

\def\nm{\noalign{\medskip}}
\newcommand{\eqnref}[1]{(\ref {#1})}

\newcommand{\Be}{\mathbf{e}}
\newcommand{\Bf}{\mathbf{f}}
\newcommand{\Bh}{\mathbf{h}}
\newcommand{\Bg}{\mathbf{g}}
\newcommand{\Bn}{{\mathbf n}}
\newcommand{\Bu}{\mathbf{u}}
\newcommand{\Bv}{\mathbf{v}}
\newcommand{\Bw}{\mathbf{w}}

\newcommand{\Bz}{\mathbf{z}}

\newcommand{\BA}{\mathbf{A}}
\newcommand{\BB}{\mathbf{B}}
\newcommand{\BE}{\mathbf{E}}
\newcommand{\BI}{\mathbf{I}}
\newcommand{\BU}{\mathbf{U}}
\newcommand{\BN}{\mathbf{N}}

\newcommand{\GG}{\mathbb{G}}
\newcommand{\CC}{\mathbb{C}}

\newcommand{\Ga}{\alpha}
\newcommand{\Gb}{\beta}
\newcommand{\Gd}{\delta}
\newcommand{\Ge}{\epsilon}
\newcommand{\Gk}{\kappa}
\newcommand{\Gl}{\lambda}
\newcommand{\Gm}{\mu}
\newcommand{\Gn}{\nu}
\newcommand{\Gs}{\sigma}
\newcommand{\GL}{\Lambda}
\newcommand{\GO}{\Omega}

\newcommand{\vk}{\Bv}
\newcommand{\vp}{\varphi}

\newcommand{\Kcal}{\mathcal{K}}
\newcommand{\Lcal}{\mathcal{L}}
\newcommand{\Scal}{\mathcal{S}}

\newcommand{\ds}{\displaystyle}
\newcommand{\la}{\langle}
\newcommand{\ra}{\rangle}

\newcommand{\p}{\partial}
\newcommand{\pd}[2]{\frac {\p #1}{\p #2}}
\newcommand{\pf}{\medskip \noindent {\sl Proof}. \ }
\newcommand{\qed}{\hfill $\Box$ \medskip}
\newcommand{\RR}{\mathbb{R}}

\newcommand{\beq}{\begin{equation}}
\newcommand{\eeq}{\end{equation}}

\DeclareMathAlphabet{\itbf}{OML}{cmm}{b}{it}

\def\by{{{\itbf y}}}
\def\bx{{{\itbf x}}}

\def\bz{{{\itbf z}}}

\def\bu{{{\itbf u}}}

\begin{document}
\title{Strong convergence of the solutions of the linear elasticity and uniformity of asymptotic expansions in the presence of small inclusions\thanks{\footnotesize This work was
supported by the ERC Advanced Grant Project MULTIMOD--267184 and
NRF grants No. 2009-0090250, 2010-0004091, and 2010-0017532.}}

\author{Habib Ammari\thanks{\footnotesize Department of Mathematics and Applications, Ecole Normale Sup\'erieure,
45 Rue d'Ulm, 75005 Paris, France (habib.ammari@ens.fr).} \and
Hyeonbae Kang\thanks{Department of Mathematics, Inha University, Incheon
402-751, Korea (hbkang@inha.ac.kr, kskim@inha.ac.kr, hdlee@inha.ac.kr).}  \and
Kyoungsun Kim\footnotemark[3] \and Hyundae Lee\footnotemark[3]}

\date{}
\maketitle

\begin{abstract}
We consider the Lam\'e system of linear elasticity when the
inclusion has the extreme elastic constants. We show that the
solutions to the Lam\'e system converge in appropriate $H^1$-norms when the
shear modulus tends to infinity (the other modulus, the compressional modulus is
fixed), and when the bulk modulus and the shear modulus tend to
zero. Using this result, we show that the asymptotic expansion of the
displacement vector in the presence of small inclusion is uniform
with respect to Lam\'e parameters.
\end{abstract}


\section{Introduction}

Recently there is growing interest in partial differential
equations with high contrast coefficients in various contexts.
Among them are the photonic and phononic band gap problems where
the electromagnetic parameter or the bulk modulus tend to
infinity, biomedical imaging where anomalous tissues have large
material parameters, and the stress concentration in between two
inclusions with extreme material properties to name a few. See
\cite{AKL09, book2, ACKLY} and references therein. The purpose of
this paper is to prove two basic theorems in relation to the PDEs
with high contrast coefficients. The first one is to show that
when the material parameters tend to the extreme, the
corresponding solutions converge strongly in appropriate norms.
The other one is to show that the asymptotic expansion of the
solution in the presence of small inclusions holds uniformly with
respect to material parameters. We prove these facts in the
context of the system of linear elasticity. Corresponding results
for the conductivity equation (a scalar equation) have been
obtained in \cite{FV89, BLY1, NV}.

We consider a linear isotropic elastic body containing an
inclusion with different elastic parameters. When the
bulk and shear moduli of the inclusion are finite, the solution
satisfies the transmission condition along the interface (the
boundary of the inclusion). If the shear modulus of the inclusion
is infinity, then the interface transmission condition is replaced
by a null condition of the displacement (see Section 2). If the
bulk and shear moduli are zero, then it is replaced by the
traction zero condition on the boundary of the inclusion. The
first objective of this paper is to prove the convergence in an
appropriate $H^1$ space of the solution to the Lam\'e system as
the bulk and shear moduli tend to the extreme (zero or infinity)
(see Theorem \ref{thm1}).

The second objective of this paper is to prove a closely related
problem of uniformity of the asymptotic expansion. In imaging
small inclusions from boundary measurements, it is of fundamental
importance to catch the boundary signature of the presence of
anomalies. In this respect, an asymptotic expansion of the
boundary perturbations of the solutions due to the presence of the
inclusion, as the diameter of the inclusion tends to zero, has
been derived. The asymptotic expansion is derived in \cite{FV89,
CMV, AK03} for the conductivity (scalar) equation, and in
\cite{AKNT02, LS02} for the Lam\'e system of linear isotropic
elasticity. The asymptotic expansions have been effectively used
for imaging diametrically small inclusions. See for example
\cite{CMV, KKL03, book1, AK06, KKL07, ABCTF}. We also mention the
topological derivative based shape optimization where the
asymptotic expansion is an essential ingredient (see for example
\cite{CGGM, GGM, AGJK-Top, ABGJKW_arxiv}). In \cite{ABGJKW_arxiv}
topological derivative based detection algorithms for the
localization of an elastic inclusion of vanishing characteristic
size have been developed and their resolution and stability with
respect to measurement and medium noises analyzed.

In these applications, it is important to know that the asymptotic
expansion holds uniformly with respect to the pair of Lam\'e
parameters. We prove this in the second half of this paper under
the assumption that the compressional modulus is bounded, which is
necessary (see Theorem \ref{firstthm} for precise statements). It
is worthwhile to mention that this result may have a relation with
the cloaking as discussed in \cite{NV, gang}.

The methods of this paper are different from those of \cite{NV}, where uniform validity of the asymptotic expansion for the conductivity (scalar) equations is proved, in that they are based on the layer potential
techniques. The solutions to the Lam\'e system can be expressed as
a single layer potential on the boundary of inclusion. We show
that $H^{-1/2}$-norms of the potentials are bounded uniformly with
respect to Lam\'e parameters, and the main results follow from
this fact.

This paper is organized as follows. In section 2, we set up the
problems for finite and extreme moduli, and review the
representation of solutions using layer potential techniques.  In
section 3, we prove that the energy functional is uniformly
bounded. As a consequence, we obtain that the potentials on the
boundary of the inclusion are uniformly bounded. In section 4 we
show that these potentials converge as the bulk and shear moduli
tend to extreme values and prove Theorem \ref{thm1}. In section 5,
we briefly discuss that similar boundedness and convergence result
hold to be true for the boundary value problem.  Section 6 is to
prove Theorem \ref{firstthm} which asserts that the small volume
expansion holds independently of Lam\'e parameters. The results
and methods hold to be true even if there are multiple inclusions.
We make a brief remark on this in the last section.

\section{Problem setting and representation of solutions}

Let $D$ be an elastic inclusion which is a bounded domain in
$\RR^d\ (d=2,3)$ with the Lipschitz boundary. Let $(\Gl, \Gm)$ be
the pair of Lam\'e (shear and compressional) parameters of $D$ while $(\Gl_0, \Gm_0)$ is
that of the background $\RR^d \setminus D$. Then the elasticity
tensors for the inclusions and the background can be written
respectively as $\CC^1= (C^1_{ijk\ell})$ and $\CC^0=
(C^0_{ijk\ell})$ where
 \begin{align*}
 C^1_{ijk\ell} &= \Gl \delta_{ij}\delta_{k\ell} + \Gm (\delta_{ik}\delta_{j\ell}+\delta_{i\ell}\delta_{jk}), \\
 C^0_{ijk\ell} &= \Gl_0 \delta_{ij}\delta_{k\ell} + \Gm_0 (\delta_{ik}\delta_{j\ell}+\delta_{i\ell}\delta_{jk}),
 \end{align*}
and the elasticity tensor for
$\RR^d$ in the presence of the inclusion $D$ is given by
\beq
 \CC:=1_{D} \CC^1 + (1-1_{D}) \CC^0,
 \eeq
where $1_D$ is the indicator function of $D$. We assume that the strong convexity condition holds, {\it
i.e.},
\begin{equation}\label{strongcon}
\Gm > 0, \quad d \Gl + 2\Gm >0,\quad \Gm_0 > 0, \quad \mbox{and} \quad d \Gl_0 + 2\Gm_0 >0\;.
\end{equation}
We also assume that
 \beq\label{minorcond}
 (\Gl - \Gl_0)(\Gm - \Gm_0) > 0,
 \eeq
which is required to have the representation of the displacement
vectors by the single layer potential in the following. We also
denote the bulk modulus by $\Gk$ which is given by $\Gk= \Gl+
2\mu/d$.

We consider the problem of the Lam\'e system of the linear
elasticity: For a given function $\Bh$ satisfying $\nabla\cdot \CC^0 \nabla^s \Bh=0$ in $\RR^d$,
\begin{equation}\label{main-eqn}
\left\{
\begin{array}{ll}
\nabla\cdot \CC \nabla^s \Bu=0 \quad & \mbox{in } \RR^d,  \\
\nm \Bu(\bx)- \Bh(\bx)=O(|\bx|^{1-d}) \quad & \mbox{as
}|\bx|\rightarrow \infty,
\end{array}
\right.
\end{equation}
where $\nabla^s \Bu$ is the symmetric gradient (or the strain tensor), {\it i.e.},
 $$
 \nabla^s \Bu:= \frac{1}{2} (\nabla \Bu + (\nabla \Bu)^{T}) \quad (T \mbox{ for
 transpose}).
 $$
Let
 $$
 \Lcal_{\Gl, \Gm} \Bu : =\nabla\cdot \CC^1 \nabla^s \Bu= \Gm \Delta \Bu + (\Gl +
\Gm) \nabla \nabla \cdot \Bu \;
 $$
and define the corresponding conormal derivative ${\partial \Bu}/{\partial \Gn}$ on $\p D$ by
 \begin{equation} \label{Lame-2}
 \pd{\Bu}{\Gn} := \CC^1 (\nabla^s \Bu) \, \Bn = \Gl (\nabla\cdot \Bu) \Bn + \Gm (\nabla
 \Bu + (\nabla \Bu)^T ) \Bn \quad \mbox{on } \p D ,
 \end{equation}
where $\Bn$ is the outward unit normal to $\p D$. Let $D^c:=\RR^d
\setminus \overline{D}$. Let $\Lcal_{ \Gl_0 , \Gm_0 }$ and
$\pd{}{\Gn_0}$ be those corresponding to $(\Gl_0 , \Gm_0)$. Then
\eqnref{main-eqn} is equivalent to the following problem:
 \beq\label{trans}
 \begin{cases}
 \ds \Lcal_{ \Gl_0 , \Gm_0 } \Bu = 0  \quad \mbox{in } D^c, \\
 \ds \Lcal_{\Gl , \Gm} \Bu = 0 \quad \mbox{in } D, \\
 \ds \Bu |_{-} = \Bu |_{+} \quad \mbox{on } \p D, \\
 \ds \pd{\Bu}{\Gn} \bigg |_{-} = \pd{\Bu}{\Gn_0} \bigg |_{+}
 \quad \mbox{on } \p D, \\
 \Bu(\bx)- \Bh(\bx)=O(|\bx|^{1 - d}) \quad \mbox{as }|\bx|\rightarrow \infty ,
 \end{cases}
 \eeq
where the subscripts $+$ and $-$ indicate the limits from outside and inside $D$, respectively.

We also consider the two limiting cases of \eqnref{trans}: when
both $\Gk$ and $\mu$ tend to $0$, and when $\mu \to \infty$ while
$\Gl$ is fixed. In relation to the latter case it is worth
mentioning that if $\Gl \to \infty$ and $\mu$ is fixed,
\eqnref{trans} approaches to a different problem. Roughly
speaking, if $\Gl \to \infty$, then $\nabla \cdot \Bu$ is
approaching to $0$ while $\Gl \nabla \cdot \Bu$ stays bounded. So
\eqnref{trans} approaches to the modified Stokes' problem $\mu
\Delta \Bu + \nabla p=0$ with $p=\Gl \nabla \cdot \Bu$. (See
\cite{AGKL}.) Hence we assume that $\Gl$ is bounded throughout
this paper.

If $\Gk=\mu=0$ (or $\Gl=\mu=0$), one can easily see what the
limiting problem should be. Since $\pd{\Bu}{\Gn} |_{-}=0$, we have
from the fourth line of \eqnref{trans} that $\pd{\Bu}{\Gn_0}
|_{+}=0$. So the elasticity equation in this case is
 \begin{equation}\label{eqnhole}
 \begin{cases}
 \ds \Lcal_{ \Gl_0 , \Gm_0 } \Bu = 0
 \quad \mbox{in } D^c, \\
 \ds \pd{\Bu}{\Gn_0} \bigg |_{+} =0
 \quad \mbox{on } \p D, \\
 \Bu(\bx)- \Bh(\bx)=O(|\bx|^{1 - d}) \quad \mbox{as }|\bx|\rightarrow \infty.
 \end{cases}
 \end{equation}

To describe the equation when $\Gm = \infty$ while $\Gl$ remains
bounded, we need to introduce the following functional space: Let
$\Psi$ be the $d(d+1)/2$ dimensional vector space defined by
\beq\label{Psi} \Psi:=\{~ \psi = (\psi^{(1)}, \cdots,
\psi^{(d)})^T ~:~ \p_i \psi^{(j)} + \p_j \psi^{(i)} = 0, \ 1\leq
i,j \leq d ~\}. \eeq We emphasize that $\Psi$ is the space of
solutions $\Lcal_{\Gl,\Gm}\Bu=0$ in $D$ and $\p\Bu/{\p \Gn}=0$ on
$\p D$ for any $(\Gl,\mu)$. Let $\psi_j$, $j=1, \ldots, d(d+1)/2$,
be a basis of $\Psi$. If $\Gm \to \infty$, then from the second
and fourth equations in \eqnref{trans} we have
 $$
 \Delta \Bu + \nabla \nabla \cdot \Bu=0 \quad \mbox{in } D, \quad
 (\nabla \Bu + (\nabla \Bu)^T ) \, \Bn =0 \quad \mbox{on } \p D,
 $$
which is another elasticity equation (with $\Gm=1$ and $\Gl=0$)
with zero traction  on the boundary. Thus there are constants
$\Ga_{j}$ such that
 $$
 \Bu(\bx)= \sum_{j=1}^{d(d+1)/2} \Ga_{j} \psi_j (\bx), \quad \bx \in D.
 $$
So, the elasticity problem when $\Gm = \infty$ is
 \begin{equation} \label{eqnhard}
 \begin{cases}
 \ds \Lcal_{ \Gl_0 , \Gm_0 } \Bu = 0 \quad \mbox{in } D^c, \\
 \ds \Bu = \sum_{j=1}^{d(d+1)/2} \Ga_{j} \psi_j  \quad \mbox{on } \p D, \\
 \Bu(\bx)- \Bh(\bx)=O(|\bx|^{1 - d}) \quad  \mbox{as }|\bx|\rightarrow \infty.
 \end{cases}
 \end{equation}
We need extra conditions to determine the coefficients $\Ga_j$.
Note that the solution $\Bu$ to \eqnref{trans} satisfies
 \beq\label{zerocon2}
 \int_{\p D} \pd{\Bu}{\Gn_0} \bigg |_{+} \cdot \psi_l \, d\Gs = 0, \quad l=1,\cdots,\frac{d(d+1)}{2}.
 \eeq
So, by taking a (formal) limit, one can expect that the solution
$\Bu$ to \eqnref{eqnhard} should satisfy the same condition, and
the constants $\Ga_{j}$ in \eqnref{eqnhard} are determined by this
orthogonality condition.

We now review the representation of the solution to \eqnref{trans} using the single layer potential for the Lam\'e system following \cite{kup, DKV, ES93, book2}. The Kelvin matrix of the fundamental solution ${\bf \Gamma} = ( \Gamma_{ij} )_{i,j=1}^d$ to the Lam{\'e} system $\Lcal_{\lambda, \mu}$ is given by
 $$
\Gamma_{ij} (\bx) : = \begin{cases}
\ds \frac{\Ga}{2 \pi} \delta_{ij} \ln |\bx| - \frac{\Gb}{2 \pi} \frac{x_i
 x_j}{|\bx|^2} & \quad \mbox{if } d=2 \;, \\
\ds -\frac{\Ga}{4\pi} \frac{\delta_{ij}}{|\bx|} - \frac{\Gb}{2 \pi} \frac{x_i
 x_j}{|\bx|^3} & \quad \mbox{if } d=3 \;,
\end{cases} \quad  \bx \neq 0\
 $$
where
 $$
 \Ga= \frac{1}{2} \left ( \frac{1}{\mu} + \frac{1}{2\mu + \lambda} \right )
 \quad\mbox{and}\quad \Gb= \frac{1}{2} \left ( \frac{1}{\mu} - \frac{1}{ 2 \mu + \lambda} \right )\;.
 $$

When $D$ is a simply connected domain, the single  layer potentials of the density function
$\vp=\begin{bmatrix} \vp_1 & \cdots & \vp_d \end{bmatrix}^T$ on $\p D$ associated with the Lam{\'e} parameters
$(\Gl, \Gm)$ are defined by
 \begin{align}
& \Scal_D [\vp] (\bx) := \int_{\p D} {\bf \Gamma} (\bx-\by)
 \vp (\by)\, d \sigma (\by)\;, \quad \bx \in \RR^d\;.
 \end{align}
We denote by $\bf\Gamma^0$, $\Scal_D^0[\vp]$  the fundamental
solution and  the single layer potential associate with the
Lam{\'e} parameter $(\Gl_0,\Gm_0)$ respectively. The conormal
derivative of $\Scal_D[\vp]$ enjoys the jump relation on $\p D$:
\beq\label{jump} \frac{\p}{\p\Gn}\Scal_D[\vp]\Big|_{\pm} = (\pm
\frac{1}{2}I + \Kcal_D^*)[\vp] \quad \mbox{a.e. on } \p D, \eeq
where $\Kcal_D^*$ is defined by
\beq
\Kcal_D^*[\vp](\bx) =
\mbox{p.v.} \int_{\p D} \frac{\p}{\p \Gn_{\bx}} {\bf \Gamma}(\bx -
\by)\vp (\by)\, d \Gs (\by) \quad \mbox{a.e. } \bx \in \p D,
\eeq
where p.v. stands for the Cauchy principal value. We denote by $\Kcal_D^{0*}$ the operator corresponding to
$(\Gl_0,\Gm_0)$.

We introduce a weighted norm, $\| \Bu \|_{H^1_w(\GO)}$, in two dimensions: let $\GO$ be either $\RR^d$ or $D^c$, and let
\begin{equation*}
\| \Bu \|^2_{H^1_w(\GO)} := \int_{\GO} \frac{|\Bu(\bx)|^2}{\sqrt{1+|\bx|^2}} d\bx + \int_{\GO} |\nabla \Bu (\bx)|^2 d\bx.
\end{equation*}
This weighted norm is introduced because the solutions $\Bu$
satisfies only $\Bu(\bx)=O(|\bx|^{-1})$ in two dimensions as
$|\bx| \to \infty$. For convenience in presenting results of this
paper, we put $W(\GO):=H^1_w(\GO)$ in two dimensions, and
$W(\GO):= H^1(\GO)$, the usual Sobolev space, in three dimensions.
Let $\Psi$ be the space introduced in \eqnref{Psi}, and define
\beq H^{-1/2}_{\Psi}(\p D) := \{~ \vp \in H^{-1/2}(\p D)^d ~:~
\la \vp , \psi \ra = 0 \mbox{ for all } \psi \in \Psi ~\}. \eeq
Here $\la \vp , \psi \ra$ denotes the $H^{-1/2}-H^{1/2}$ product.
Then $\pm \frac{1}{2}I + \Kcal_D^*$ is invertible on
$H^{-1/2}_{\Psi}(\p D)$. We also have \beq\label{Scalest} \|
\Scal_D[\vp] \|_{W(\RR^d)} \le C \| \vp \|_{H^{-1/2}(\p D)} \eeq
for all $\vp \in H^{-1/2}(\p D)^d$.

It is proved in \cite{ES93} that the solution $\Bu$ to \eqnref{trans} is represented as
\beq\label{reprek}
\Bu(\bx) =\begin{cases}
\Bh(\bx) + \Scal_D^0[\vp](\bx),  \quad &\bx \in D^c , \\
\Scal_D [\psi](\bx), \quad &\bx \in D,
\end{cases}
\eeq
where the pair $(\vp,\psi) \in H_{\Psi}^{-1/2}(\p D)\times H^{-1/2}(\p D) $ is the solutions to
\beq \label{int_eq}
\begin{cases}
\ds \Scal_D[\psi](\bx) - \Scal_D^0[\vp](\bx) = \Bh(\bx) \\
\nm
\ds \pd{\Scal_D[\psi]}{\Gn}\Big|_{-}(\bx) - \pd{\Scal_D^0 [\vp]}{\Gn_0}\Big|_{+}(\bx) = \pd{\Bh}{\Gn_0}(\bx)
\end{cases}
\mbox{for } \bx \in \p D.
\eeq

  Even if $\Gk=\mu=0$ or $\mu=\infty$, we have a similar representation:
\beq\label{repreinfty}
\Bu_\Gk(\bx) = \Bh(\bx) + \Scal_D^0[\vp_\Gk](\bx),  \quad \bx \in D^c, \ \Gk=0, \infty.
\eeq
When $\Gk=\mu=0$,  $\vp_0$ satisfies
\beq\label{inteqzero}
\left(\frac{1}{2}I + (\Kcal_D^{0})^*\right)[\vp_0] = -\pd{\Bh}{\Gn_0} \quad\mbox{on } \p D,
\eeq
and if $\mu=\infty$, then $\vp_\infty$ satisfies
\beq\label{inteqinf}
\left(-\frac{1}{2}I + (\Kcal_D^{0})^* \right)[\vp_\infty] = -\pd{\Bh}{\Gn_0} \quad\mbox{on } \p D.
\eeq
We emphasize that $\vp_{\Gk} \in H^{-1/2}_{\Psi}(\p D)$. See, for example, \cite{book2} for details of the above mentioned representation of the solutions.

A similar representation formula holds for the solutions to the
boundary value problems. Let $\GO$ be a bounded Lipschitz domain
in $\RR^d$ containing $D$, which is also Lipschitz. Let $\Bu$ be
the solution to \beq\label{boundaryvalue} \nabla\cdot \CC \nabla^s
\Bu=0 \quad \mbox{in } \GO, \eeq with either the Dirichlet
boundary condition $\Bu=\Bf$ or the Neumann boundary condition
$\pd{\Bu}{\nu_0}=\Bg$ on $\p\GO$. Let \beq\label{defBh} \Bh(\bx):=
- \int_{\p\GO} {\bf \Gamma^0} (\bx-\by) \pd{\Bu}{\nu_0} \Big|_{-}
(\by) \, d \sigma (\by) + \int_{\p\GO} \pd{{\bf \Gamma^0}
(\bx-\by)}{\nu_0(\by)} \Bu (\by) \, d \sigma (\by), \quad \bx \in
\GO. \eeq Then $\Bu$ is represented as \beq\label{reprek-bvp}
\Bu(\bx) =\begin{cases}
\Bh(\bx) + \Scal_D^0[\vp](\bx),  \quad &\bx \in \Omega\setminus \overline{D}, \\
\Scal_D [\psi](\bx), \quad &\bx \in D,
\end{cases}
\eeq
where  the pair $(\vp,\psi) \in H_{\Psi}^{-1/2}(\p D)\times H^{-1/2}(\p D) $ is the solutions to \eqnref{int_eq}.

\section{Energy estimates}

Let
\beq\label{funcJ}
J[\Bu]:= \frac{1}{2}\int_{D} \CC^1 \nabla^s \Bu:\nabla^s \Bu
+ \frac{1}{2}\int_{D^c} \CC^0 \nabla^s(\Bu - \Bh):\nabla^s(\Bu - \Bh).
\eeq
Here and throughout this paper $\BA:\BB= \sum_{i,j=1}^d a_{ij} b_{ij}$ for $\BA=(a_{ij})$ and $\BB=(b_{ij})$. For the solution $\Bu$ to \eqnref{trans}, we prove that $J[\Bu]$ is bounded regardless of $\Gk$ and $\mu$. More precisely we prove the following lemma.
\begin{lem}\label{lem-infty}
Let $\Bu$ be the solution to \eqnref{trans}. If $\Gl \le \GL$ for
some constant $\GL$, then there is a constant $C$ depending on
$\GL$, but otherwise independent of $\mu$ and $\Gk$, such that
\beq\label{JGkuk} J[\Bu] \leq C \left\| \pd{\Bh}{\Gn_0}
\right\|_{H^{-1/2}(\p D)}^2. \eeq
\end{lem}

As a consequence of Lemma \ref{lem-infty} we have
\begin{lem}\label{potential}
Let $\vp$ be the potential defined in \eqnref{reprek}. If $\Gl \le
\GL$ for some constant $\GL$, then there is a constant $C$
depending on $\GL$, but otherwise independent of $\mu$ and $\Gk$,
such that \beq\label{JGkuk2} \| \vp \|_{H^{-1/2}(\p D)} \leq C
\left\| \pd{\Bh}{\Gn_0} \right\|_{H^{-1/2}(\p D)}. \eeq
\end{lem}
\pf Let $\Bv:= \Bu -\Bh$. Then \eqnref{reprek} yields  $\Bv(\bx) =
\Scal_D^0[\vp](\bx)$ for $\bx \in D^c$. Thus, we have from
\eqnref{jump}
$$
\pd{\Bv}{\Gn_0}\Big|_{+} = \Big( \frac{1}{2}I + (\Kcal_D^0)^*\Big)[\vp] \quad \mbox{on } \p D.
$$
Since $\frac{1}{2}I + (\Kcal_D^0)^*$ is invertible on $H^{-1/2}_\Psi(\p D)$, we have
$$
\| \vp \|_{H^{-1/2}(\p D)} \leq C \left\| \pd{\Bv}{\Gn_0}\Big|_{+} \right\|_{H^{-1/2}(\p D)}.
$$

Let $\eta$ be a function in $H^{1/2}(\p D)$ satisfying $\int_{\p
D} \eta=0$ and let $\Bw$ be the solution to $\Delta \Bw=0$ in
$D^c$ with $\Bw(\bx)=O(|\bx|^{1-d})$ and $\Bw=\eta$ on $\p D$, so
that the following estimate holds:
$$
\| \nabla^s \Bw \|_{L^2(D^c)} \le C \| \eta \|_{H^{1/2}(\p D)}.
$$
Since
\beq
\int_{D^c} \CC^0 \nabla^s \Bv: \nabla^s \Bw d\bx = -\int_{\p D} \pd{\Bv}{\Gn_0}\Big|_{+} \cdot \eta d\Gs(\bx),
\eeq
we have
\begin{align*}
\Big| \int_{\p D}  \pd{\Bv}{\Gn_0}\Big|_{+} \cdot \eta d\Gs(\bx) \Big|
& \leq C \|\nabla^s \Bv\|_{L^2(D^c)} \| \nabla^s \Bw \|_{L^2(D^c)}\\
& \leq C \|\nabla^s \Bv\|_{L^2(D^c)} \| \eta \|_{H^{1/2}(\p D)}.
\end{align*}
Since $\eta \in H^{1/2}(\p D)$ is arbitrary, we have from \eqnref{JGkuk}
\begin{equation*}
\Big\| \pd{\Bv}{\Gn_0}\Big|_{+} \Big\|_{H^{-1/2}(\p D)} \leq C \|\nabla^s \Bv\|_{L^2(D^c)} \leq C \left\| \pd{\Bh}{\Gn_0} \right\|_{H^{-1/2}(\p D)},
\end{equation*}
and so follows \eqnref{JGkuk2}. \qed

\noindent{\sl Proof of Lemma \ref{lem-infty}}.
Let $\vk:= \Bu -\Bh$. It is known (see for example \cite{CK}) that $\vk$ is the minimizer in $W(\RR^d)$ of the functional
\beq\label{IGk}
I[\Bv]:=
\frac{1}{2}\int_{\RR^d} \CC (\nabla^s \Bv + 1_{D} \GG \nabla^s \Bh) : (\nabla^s \Bv + 1_{D} \GG \nabla^s \Bh),
\eeq
where $\GG=\BI_4 - (\CC^1)^{-1} \CC^0$ and $\BI_4$ is the identity 4-tensor. Note that
\begin{align}
I[\Bv] & =\frac{1}{2}\int_{D} (\CC^1 \nabla^s(\Bv + \Bh) - \CC^0\nabla^s \Bh):(\nabla^s(\Bv + \Bh) - (\CC^1)^{-1} \CC^0\nabla^s \Bh) \nonumber \\
& \qquad + \frac{1}{2}\int_{D^c}\CC^0\nabla^s \Bv:\nabla^s \Bv \nonumber \\
& = \frac{1}{2}\int_{D} \CC^1 \nabla^s \Bv: \nabla^s \Bv
+ \int_{D}(\CC^1 - \CC^0) \nabla^s \Bv: \nabla^s \Bh \nonumber \\
& \qquad + \frac{1}{2}\int_{D} (\CC^1 - 2\CC^0 + \CC^0 (\CC^1)^{-1} \CC^0 )\nabla^s \Bh:\nabla^s \Bh  + \frac{1}{2}\int_{D^c} \CC^0 \nabla^s \Bv :\nabla^s \Bv. \label{IGk2}
\end{align}

Let
\beq
\Bv_{\infty} := \Bu_{\infty} - \Bh.
\eeq
Then $\Bv_{\infty} \in W(\RR^d)$, and $(\Bv_{\infty} + \Bh)|_{D} \in \Psi$ which implies that $\nabla^s (\Bv_{\infty} + \Bh)=0$ in $D$. So, we have from the first line in \eqnref{IGk2} that
\beq\label{Ivinfty}
I[\Bv_{\infty}] = \frac{1}{2}\int_{D} \CC^0 (\CC^1)^{-1} \CC^0 \nabla^s \Bh:\nabla^s \Bh + \frac{1}{2}\int_{D^c} \CC^0 \nabla^s \Bv_{\infty} : \nabla^s \Bv_{\infty}.
\eeq

We then have
\begin{align*}
J [\Bu] &= J[\vk +\Bh] \\
& = \frac{1}{2}\int_{D}\CC^1 \nabla^s \vk:\nabla^s \vk + \int_{D}\CC^1 \nabla^s \vk:\nabla^s \Bh + \frac{1}{2}\int_{D}\CC^1 \nabla^s \Bh:\nabla^s \Bh\\
& \quad + \frac{1}{2}\int_{D^c} \CC^0 \nabla^s \vk:\nabla^s \vk\\
& = I[\vk] + \int_{D} \CC^0 \nabla^s \vk:\nabla^s \Bh + \int_{D} \CC^0 \nabla^s \Bh:\nabla^s \Bh - \frac{1}{2}\int_{D} \CC^0 (\CC^1)^{-1} \CC^0 \nabla^s \Bh:\nabla^s \Bh \\
& = I[\vk] + \int_{D} \CC^0 \nabla^s \Bu:\nabla^s \Bh - \frac{1}{2}\int_{D} \CC^0 (\CC^1)^{-1} \CC^0 \nabla^s \Bh:\nabla^s \Bh.
\end{align*}
Since $I[\vk] \le I[\Bv_{\infty}]$, it follows from \eqnref{Ivinfty} that
\begin{align}
J [\Bu] & \leq I[\Bv_{\infty}] + \int_{D} \CC^0 \nabla^s \Bu:\nabla^s \Bh - \frac{1}{2}\int_{D} \CC^0 (\CC^1)^{-1} \CC^0 \nabla^s \Bh:\nabla^s \Bh \nonumber \\
& = \frac{1}{2}\int_{D^c} \CC^0 \nabla^s \Bv_{\infty}:\nabla^s \Bv_{\infty} + \int_{D} \CC^0 \nabla^s \Bu:\nabla^s \Bh . \label{JGk}
\end{align}
Note that since $\Gl$ is bounded, we have
\begin{align*}
\int_{D} \CC^0 \nabla^s \Bu:\nabla^s \Bh &= C( \int_{D} |\nabla^s \Bu |^2 + \int_D |\nabla^s \Bh|^2 ) \\
 & \le C \left( \frac{1}{\mu} \int_{D} \CC^1 \nabla^s \Bu:\nabla^s \Bu + \int_D |\nabla^s \Bh|^2 \right) \\
 & \le C \left( \frac{1}{\mu} J [\Bu] + \int_D |\nabla^s \Bh|^2 \right).
\end{align*}
So, if $\mu$ is sufficiently large, then we have from \eqnref{JGk} that
\begin{equation*}
J[\Bu] \leq C\Big(\|\nabla^s \Bv_{\infty}\|^2_{L^2(D^c)} + \|\nabla^s \Bh \|^2_{L^2(D)}\Big)
\end{equation*}
for some constant $C$. Since
$$
\|\nabla^s \Bv_{\infty}\|_{L^2(D^c)} \le C \| \vp_\infty \|_{H^{-1/2}(\p D)} \le C \left\| \pd{\Bh}{\Gn_0} \right\|_{H^{-1/2}(\p D)},
$$
we have \eqnref{JGkuk} when $\mu$ is large.

When $\Gk$ and $\mu$ are bounded, we need a function which plays the role of $\Bv_\infty$ in the above. For that we use
$\vp_0$ in \eqnref{repreinfty}: define
 \beq\label{repzero}
 \Bv_0 (\bx) := \Scal_D^0 [\vp_0](\bx) \quad \mbox{for }\bx \in \RR^d.
 \eeq
It is worth emphasizing that $\Bv_0$ is defined not only on $D^c$
but on $\RR^d$. Then one can show as above that
$$
J[\Bu]
\leq I[\Bv_0] + \int_{D} \CC^0 \nabla^s \vk:\nabla^s \Bh + \int_{D} \CC^0 \nabla^s \Bh:\nabla^s \Bh - \frac{1}{2}\int_{D} \CC^0 (\CC^1)^{-1} \CC^0 \nabla^s \Bh:\nabla^s \Bh.
$$
Using \eqnref{IGk2}, one can see that
\begin{align}
J[\Bu]
& \le \frac{1}{2}\int_D \CC^1 \nabla^s \Bv_0 : \nabla^s \Bv_0 + \int_D (\CC^1 -\CC^0) \nabla^s \Bv_0 : \nabla^s \Bh + \frac{1}{2}\int_D \CC^1 \nabla^s \Bh : \nabla^s \Bh \nonumber\\
& \quad + \frac{1}{2}\int_{D^c} \CC^0 \nabla^s \Bv_0 : \nabla^s \Bv_0 + \int_D \CC^0 \nabla^s \vk : \nabla^s \Bh. \label{JGk2}
\end{align}
Since $\frac{\p (\Bv_0 + \Bh)}{\p \Gn_0}\Big|_{+} = 0$ on $\p B$, we have
\begin{align*}
\int_D \CC^0 \nabla^s \vk : \nabla^s \Bh & = \int_{\p D} \pd{\Bh}{\Gn_0} \cdot \vk = -\int_{\p D} \pd{\Bv_0}{\Gn_0}\Big|_{+} \cdot \vk
\\
& = \int_{D^c} \CC^0 \nabla^s \Bv_0 : \nabla^s \vk\\
& \leq \frac{C}{\Ge} \int_{D^c} \CC^0 \nabla^s \Bv_0 : \nabla^s \Bv_0 + C\Ge \int_{D^c} \CC^0 \nabla^s \vk : \nabla^s \vk
\end{align*}
for a (small) constant $\Ge$. If $\Ge$ is sufficiently small, then we obtain by combining this with \eqnref{JGk2}
\begin{equation}\label{bound-zero}
J[\Bu] \leq C (\|\nabla^s \Bv_0 \|_{L^2(\RR^d)}^2 + \|\nabla^s \Bh \|_{L^2(D)}^2)
\end{equation}
for some constant $C$ independent of $\Gk$ and $\mu$. Since
$$
\|\nabla^s \Bv_{0}\|_{L^2(D^c)} \le C \| \vp_0 \|_{H^{-1/2}(\p D)} \le C \left\| \pd{\Bh}{\Gn_0} \right\|_{H^{-1/2}(\p D)},
$$
we have \eqnref{JGkuk}. This completes the proof.
\qed

\section{Convergence of potentials and solutions}\label{sec2}

Lemma \ref{potential} shows that the potential defined in
\eqnref{reprek} is uniformly bounded with respect to $\mu$ and
$\Gl$ as long as $\Gl$ is bounded. We now prove the following
lemma.
\begin{lem}\label{lem2}
Let $\vp$, $\vp_0$ and $\vp_\infty$ be potentials defined by \eqnref{reprek}, \eqnref{inteqzero}
and \eqnref{inteqinf}, respectively.
\begin{itemize}
\item[{\textrm (i)}] Suppose that $\Gl \le \GL$ for some constant $\GL$. There are constants $\Gm_1$ and $C$ such that
\begin{equation}\label{conv-vpinf}
\| \vp - \vp_\infty \|_{H^{-1/2}(\p D)} \le \frac{C}{\sqrt{\Gm}} \left\| \pd{\Bh}{\Gn_0} \right\|_{H^{-1/2}(\p D)}
\end{equation}
for all $\Gm \ge \Gm_1$.
\item[{\textrm (ii)}] There are constants $\Gd$ and $C$ such that
\begin{equation}\label{conv-vpzero}
\| \vp - \vp_0 \|_{H^{-1/2}(\p D)} \le C(\Gk+\mu)^{1/4} \left\| \pd{\Bh}{\Gn_0} \right\|_{H^{-1/2}(\p D)}
\end{equation}
for all $\Gk,\mu \le \Gd$.
\end{itemize}
\end{lem}

\pf We may assume that $\| \pd{\Bh}{\Gn_0} \|_{H^{-1/2}(\p D)} =1$.
Let $\Bw = \Bu - \Bu_{\infty}$ so that $\Bw$ satisfies
\beq
\begin{cases}
\nabla \cdot \CC^0 \nabla^s \Bw = 0 \quad \mbox{in } D^c,\\
\Bw - \Bu |_{\overline D} \in \Psi\quad \mbox{in } \overline{D},\\
\Bw(\bx) = O(|\bx|^{1-d}) \quad \mbox{as } |\bx| \rightarrow \infty.
\end{cases}
\eeq
Since $\nabla^s (\Bw - \Bu)=0$ in $D$, we have from Lemma \ref{lem-infty} that
\begin{align}
\int_D | \nabla^s \Bw |^2 d\bx & \leq \frac{1}{2\Gm} \int_D \CC^1 \nabla^s \Bw : \nabla^s \Bw d\bx \nonumber \\
& = \frac{1}{2\Gm} \int_D \CC^1 \nabla^s \Bu : \nabla^s \Bu d\bx \leq \frac{C}{\Gm}. \label{oneovermu}
\end{align}

Let $\psi_j$, $j=1, \ldots, d(d+1)/2$, be a basis of $\Psi$ as before, and let
$$
\Be=\sum_{j=1}^{d(d+1)/2}\Gb_j \psi_j,
$$
where $\Gb_j$ are chosen so that
\beq
\int_D \Big[(\Bw - \Be) \cdot \psi_j + \nabla (\Bw - \Be) : \nabla \psi_j \Big] = 0,  \quad j=1, \ldots, d(d+1)/2.
\eeq
We then apply Korn's inequality (see for example \cite{cialet}) to $\Bw - \Be$ to have
\beq
\int_D\Big[|\Bw - \Be|^2 + |\nabla (\Bw - \Be)|^2\Big]
\leq C \int_D |\nabla^s (\Bw - \Be)|^2 \leq C \int_D |\nabla^s \Bw|^2
\eeq
for some constant $C$ independent of $\Gm$. It then follows from \eqnref{oneovermu} that
\begin{equation}\label{inside-infty}
\| \Bw - \Be \|_{H^1(D)} \leq \frac{C}{\sqrt{\Gm}},
\end{equation}
and from the trace theorem that
\beq\label{onehalf}
\| \Bw - \Be \|_{H^{1/2}(\p D)} \leq \frac{C}{\sqrt{\Gm}}
\eeq
for some constant $C$ independent of $\Gm$. By the strong convexity \eqnref{strongcon} of $\CC^0$, there
is a constant $C$ such that
\begin{align*}
\|\nabla^s \Bw\|_{L^2(D^c)}^2
& \leq C \int_{D^c} \CC^0 \nabla^s \Bw : \nabla^s \Bw d\bx\\
& = - C \int_{\p D} \pd{\Bw}{\Gn_0} \Big|_{+} \cdot \Bw \ d\Gs(\bx)\\
& = - C \int_{\p D} \pd{\Bw}{\Gn_0} \Big|_{+} \cdot (\Bw - \Be) \ d\Gs(\bx)\\
& \leq C \| \Bw - \Be \|_{H^{1/2}(\p D)} \Big\| \pd{\Bw}{\Gn_0}\Big|_{+} \Big\|_{H^{-1/2}(\p D)},
\end{align*}
where the second equality holds because of the orthogonality
property \eqnref{zerocon2}. It then follows from \eqnref{onehalf}
that
\beq\label{outside-1}
\|\nabla^s \Bw\|_{L^2(D^c)}^2 \leq
\frac{C}{\sqrt{\Gm}} \Big\| \pd{\Bw}{\Gn_0}\Big|_{+}
\Big\|_{H^{-1/2}(\p D)}.
\eeq
We then obtain using the $H^{-1/2}$-$H^{1/2}$ duality, divergence theorem on $D^c$, Korn's inequality, and the trace theorem that
\begin{equation*}
\Big\| \pd{\Bw}{\Gn_0}\Big|_{+} \Big\|^2_{H^{-1/2}(\p D)} \leq C \|\nabla^s \Bw\|^2_{L^2(D^c)} ,
\end{equation*}
and from \eqnref{outside-1} that
\beq\label{412}
\Big\| \pd{\Bw}{\Gn_0}\Big|_{+} \Big\|_{H^{-1/2}(\p D)} \leq \frac{C}{\sqrt{\Gm}} .
\eeq

Using the representations \eqnref{reprek} and \eqnref{repreinfty}, we have
\beq
\Bw(\bx) = \Bu(\bx) - \Bu_{\infty}(\bx) = \Scal_D^0[\,\varphi -  \varphi_{\infty}](\bx), \quad \bx \in D^c.
\eeq
Thus, \eqnref{jump} yields
\begin{equation}\label{conormal}
\pd{\Bw}{\Gn_0}\Big|_{+} = \Big( \frac{1}{2}I + (\Kcal_D^0)^*\Big)[\,\varphi - \varphi_{\infty}] \quad \mbox{on } \p D.
\end{equation}
So, \eqnref{conv-vpinf} follows from \eqnref{412}.

To prove \eqnref{conv-vpzero}, let $\Bv_0$ be as defined in
\eqnref{repzero} and let $\Bz:= \vk - \Bv_0$ in $\RR^d$. Then
$\Bz=\Bu - \Bu_0$ in $D^c$ and the following holds
\begin{align*}
\int_{D^c} |\nabla^s \Bz|^2 d \bx
& \leq C \int_{D^c} \CC^0 \nabla^s \Bz : \nabla^s \Bz d \bx\\
& = - C \int_{\p D} \pd{\Bz}{\Gn_0}\Big|_{+} \cdot \Bz = - C \int_{\p D} \pd{\Bu}{\Gn_0}\Big|_{+} \cdot \Bz \\
& = - C \Big(\int_{\p D} \pd{\Bu}{\Gn_0}\Big|_{+} \cdot \Bu - \int_{\p D} \pd{\Bu}{\Gn_0}\Big|_{+} \cdot \Bu_0\Big)
\end{align*}
for some constant $C>0$. Since
\begin{equation*}
\int_{\p D} \pd{\Bu}{\Gn_0}\Big|_{+} \cdot \Bu = \int_{\p D} \pd{\Bu}{\Gn}\Big|_{-} \cdot \Bu = \int_D \CC^1 \nabla^s \Bu : \nabla^s \Bu \geq 0,
\end{equation*}
we have
\begin{align*}
\int_{D^c} |\nabla^s \Bz|^2 d \bx
& \leq C \int_{\p D} \pd{\Bu}{\Gn_0}\Big|_{+} \cdot \Bu_0\\
& = C \int_{\p D} \pd{\Bu}{\Gn_0}\Big|_{+} \cdot (\Bv_0 + \Bh)\\
& = C \int_{\p D} \pd{\Bu}{\Gn}\Big|_{-} \cdot (\Bv_0 + \Bh)\\
& = C \int_D \CC^1 \nabla^s \Bu : \nabla^s (\Bv_0 + \Bh).
\end{align*}
By Cauchy's inequality, we obtain that
\begin{equation*}
\int_D \CC^1 \nabla^s \Bu : \nabla^s (\Bv_0 + \Bh)
\leq C \Big(\int_D \CC^1 \nabla^s \Bu : \nabla^s \Bu\Big)^{1/2}\Big(\int_D \CC^1 \nabla^s (\Bv_0 + \Bh) : \nabla^s (\Bv_0 + \Bh)\Big)^{1/2}.
\end{equation*}
Thus, from \eqnref{JGkuk} it follows that
\begin{align*}
\int_{D^c} |\nabla^s \Bz|^2 d \bx & \leq C \Big(\int_D \CC^1 \nabla^s (\Bv_0 + \Bh) : \nabla^s (\Bv_0 + \Bh)\Big)^{1/2} \\
& \le C \sqrt{\Gk+\mu} \| \nabla^s (\Bv_0 + \Bh) \|_{L^2(D)} \le C \sqrt{\Gk+\mu}
\end{align*}
for a constant $C$ independent of $\Gk$. Therefore, we arrive at
\beq
 \Big\| \pd{\Bz}{\Gn_0}\Big|_{+} \Big\|_{H^{-1/2}(\p D)} \leq C \|\nabla^s \Bz\|_{L^2(D^c)} \leq C (\Gk+\mu)^{1/4}.
\eeq
Note that $\Bz = \Bu - \Bu_0 = \Scal_D^0[\vp-\vp_0]$ in $D^c$. So by the same reasoning as above we have \eqnref{conv-vpzero}, and the proof is complete. \qed

As a consequence of Lemma \ref{lem2}, we obtain the first main result of this paper.
\begin{thm}\label{thm1}
Suppose that \eqnref{strongcon} and \eqnref{minorcond} hold. Let
$\Bu$, $\Bu_{\infty}$ and $\Bu_0$ be the solutions to
\eqnref{trans}, \eqnref{eqnhard} and \eqnref{eqnhole},
respectively.
\begin{itemize}
\item[{\textrm (i)}] Suppose that $\Gl \le \GL$ for some constant $\GL$. There are constants $\mu_1$ and $C$ such that
\begin{equation}\label{convergence-infty}
\| \Bu - \Bu_\infty \|_{W(\RR^d)} \le \frac{C}{\sqrt{\Gm}} \left\| \pd{\Bh}{\Gn_0} \right\|_{H^{-1/2}(\p D)}
\end{equation}
for all $\mu \ge \mu_1$.
\item[{\textrm (ii)}] There are constants $\Gd$ and $C$ such that
\begin{equation}\label{convergence-zero}
\| \Bu - \Bu_0 \|_{W(D^c)} \le C(\Gk+\mu)^{1/4} \left\| \pd{\Bh}{\Gn_0} \right\|_{H^{-1/2}(\p D)}
\end{equation}
for all $\Gk, \mu \le \Gd$.
\end{itemize}
\end{thm}

It is not clear if the convergence rate, $\mu^{-1/2}$ and $(\Gk+\mu)^{1/4}$, are optimal or not.

\medskip

\noindent{\sl Proof of Theorem \ref{thm1}}.
Assume that $\| \pd{\Bh}{\Gn_0} \|_{H^{-1/2}(\p D)}=1$. Since $\Bu-\Bu_0=\Scal_D^0[\,\vp - \vp_0]$ on $D^c$, \eqnref{convergence-zero} follows from \eqnref{conv-vpzero}.

Since $\Bu-\Bu_\infty=\Scal_D^0[\,\vp - \vp_\infty]$ on $D^c$, we have
\begin{equation}\label{outside-infty}
\| \Bu-\Bu_\infty \|_{W(D^c)} \leq \frac{C}{\sqrt\Gm}.
\end{equation}
Moreover, we have
$$
\| \Bu-\Bu_\infty \|_{H^{1/2}(\p D)} = \| \Scal_D^0[\varphi -  \vp_{\infty}] \|_{H^{1/2}(\p D)} \leq C\|\varphi - \vp_{\infty}\|_{H^{-1/2}(\p D)} \leq \frac{C}{\sqrt\Gm},
$$
and hence
\beq\label{inside-infty2}
\| \Bu-\Bu_\infty \|_{H^{1}(D)} \leq \frac{C}{\sqrt\Gm}.
\eeq
This completes the proof. \qed

\section{Boundary value problems}

We now show that the results on the boundedness of the energy
functional and on the convergence of solutions similar to the
previous ones hold for the boundary value problems.

Let $\GO$ be a bounded domain in $\RR^d$ and let $D$ be an open
subset in $\GO$. We assume that $\GO$ and $D$ have Lipschitz
boundaries and satisfy
\begin{equation*}
\mbox{dist}(D, \p\GO) \geq c
\end{equation*}
for some $c >0$. We consider
\beq\label{NeumannBvp}
\begin{cases}
\nabla\cdot\CC \nabla^s \Bu=0&\quad \mbox{in  } \GO,\\
\ds\pd{\Bu}{\nu_0}=\mathbf{g} &\quad\mbox{on } \p \GO,\\
\Bu|_{\p \GO} \in L^2_\Psi ( \p \GO),
\end{cases}
\eeq where $\CC$ is given by \beq \CC = 1_D \CC^1 + 1_{\GO
\setminus \overline D} \CC^0. \eeq The Dirichlet problem can be
treated in the exactly same way.

The relevant energy functional for the boundary value problem is
\beq
J_{\GO}[\Bu]:=\frac{1}{2}\int_{\GO}\CC \nabla^s \Bu : \nabla^s \Bu.
\eeq
Then the solution $\Bu$ to \eqnref{NeumannBvp} is the minimizer of $J_{\GO}$ over $H^1(\GO)$ with the given boundary condition. Let $\Bu_\infty$ be the solution when $\mu=\infty$ ($\Gl$ is bounded). Then we have
$$
J_{\GO}[\Bu] \le J_{\GO}[\Bu_\infty] = \frac{1}{2}\int_{\GO}\CC \nabla^s \Bu_\infty : \nabla^s \Bu_\infty.
$$
Since $\Bu_\infty |_D \in \Psi$, we have $\CC^1 \nabla^s \Bu_\infty=0$ in $D$, and so,
\beq\label{uinfty}
J_{\GO}[\Bu] \le \frac{1}{2}\int_{\GO\setminus D}\CC^0 \nabla^s \Bu_\infty : \nabla^s \Bu_\infty \le C \left\|\Bg \right\|_{H^{-1/2}(\p \GO)},
\eeq
where $C$ is independent of $\mu, \Gl$.

Using \eqnref{uinfty} one can show as before that
\beq\label{conv-bvp1}
\| \Bu - \Bu_\infty \|_{H^1(\GO)} \le \frac{C}{\sqrt{\Gm}}
\eeq
for all $\mu \ge \mu_1$ when $\Gl$ is bounded, and
\beq\label{conv-bvp2}
\| \Bu - \Bu_0 \|_{H^1(\GO \setminus D)} \le C(\Gk+\mu)^{1/4}
\eeq
for all $\Gk, \mu \le \Gd$. Here $\Bu_0$ is the solution when $\Gk=\mu=0$.
We also note that \eqnref{uinfty} together with Korn's inequality implies that $\| \Bu \|_{H^{1/2}(\p\GO)} \le C$ independently of $\mu, \Gl$. It then follows from \eqnref{defBh} that
\beq\label{hbound}
\| \Bh \|_{H^{1}(\GO)} \le C
\eeq
independently of $\mu, \Gl$.

\section{Uniformity of asymptotic expansions}\label{section5}

Let us first recall the notion of elastic moment tensors (EMTs)
associated to the inclusion $D$ with the Lam\'e constants $(\Gm,
\Gl)$ when the background Lam\'e constants are $(\Gl_0, \Gm_0)$.
Let $\Ga, \beta\in \mathbb{N}^d$ be multi-indices, and let
$\{\Be_k\}_{k=1}^d$ be the standard basis of $\RR^d$. For $\Ga\in
\mathbb{N}^d,~j=1,\ldots,d$, let $(\vp_\Ga^j, \psi_\Ga^j)$ be the
solution to \eqnref{int_eq} with $\Bh$ replaced by  $\bx^\Ga
\Be_j$. Here $\bx^\Ga=x_1^{\Ga_1} \cdots x_d^{\Ga_d}$. The EMT
associated with $D$ is defined by \beq M_{\Ga\beta}^j(D) =
\int_{\p D} \bx^\beta   \varphi_\Ga^j d\sigma \eeq for
$\Ga,\beta\in\mathbb{N}^d$ and $j=1, \ldots,d$. It is worth
mentioning that $M_{\Ga\Gb}^j$ is a vector, {\it i.e.},
$M_{\Ga\Gb}^j= (m_{\Ga\Gb}^{j1}, m_{\Ga\Gb}^{j2},
m_{\Ga\Gb}^{j3})$ where \beq m_{\Ga\beta}^{jp}(D) = \int_{\p D}
\bx^\beta \Be_p \cdot  \varphi_\Ga^j d\sigma. \eeq If
$|\Ga|=|\Gb|=1$, we may write $m_{\Ga\beta}^{jp}$ as $m_{ijpq}$
for $i,j,p,q=1,\ldots, d$. It is known that $(m_{ijpq})$ is an
(anisotropic) elasticity tensor. See \cite{AKNT02}. We emphasize
that \beq\label{vpGaj} \| \vp_\Ga^j \|_{H^{-1/2}(\p D)} \leq C
\eeq for some $C$ independent of $\mu$ and $\Gk$ as long as $\Gl$
is bounded.

Suppose that $D$ is diametrically small and it is given by \beq D=
\bz_0+ \epsilon B, \eeq where $\Ge$ is a small parameter
representing the diameter of $D$, $B$ is a reference domain
containing $0$, and $\bz_0$ represents the location of $D$. Let
$\Bu$ be the solution to \eqnref{NeumannBvp} and $\BU$ be the
background solution, {\it i.e.}, the solution to \beq
\begin{cases}
\nabla\cdot\CC^0 \nabla^s \BU=0&\quad \mbox{in  } \GO,\\
\ds\pd{\BU}{\nu_0}=\mathbf{g} &\quad\mbox{on } \p \GO,\\
\BU|_{\p \GO} \in L^2_\Psi ( \p \GO).
\end{cases}
\eeq
Let $\BN$ be the Neumann function which is the solution to
\beq
\begin{cases}
\nabla\cdot\CC^0 \nabla^s \BN(\bx,\by)=-\delta_\by(\bx) \BI &\quad \mbox{in  } \GO,\\
\ds\pd{\BN(\,\cdot, \by)}{\nu_0}= - \frac{1}{|\p \GO|} \BI &\quad\mbox{on } \p \GO,\\
\BN(\,\cdot, \by)\in L^2_\Psi ( \p \GO) \quad \mbox{for each } \by\in \GO,
\end{cases}
\eeq
where $\BI$ is the identity matrix. It is proved in \cite{AKNT02} (see also \cite{book1}) that
the following asymptotic expansion holds on $\p\GO$: for $\bx \in \p\GO$
$$
\Bu(\bx)= \BU(\bx)-  \sum_{j=1}^d \sum_{|\Ga|=1}^d  \sum_{|\beta|=1}^{d+1-|\Ga|} \frac{\Ge^{|\Ga|+|\beta|+d-2}}{\Ga !\beta !} \p^\Ga U_j (\bz_0 ) \p_z^\beta \BN (\bx, \bz_0 ) M_{\Ga\beta}^j(B)
 + O(\Ge^{2d}),
$$
where $\BU=(U_1, \ldots, U_d)$. Our goal in this section is to
show that this asymptotic formula holds uniformly in $\lambda$ and
$\mu$. More precisely we have the following result.
\begin{thm} \label{firstthm}
Suppose that \eqnref{strongcon} and \eqnref{minorcond} hold. We have for $\bx \in \p\GO$
\beq \label{asympunif}
\Bu(\bx)= \BU(\bx)-  \sum_{j=1}^d \sum_{|\Ga|=1}^d  \sum_{|\beta|=1}^{d+1-|\Ga|} \frac{\Ge^{|\Ga|+|\beta|+d-2}}{\Ga !\beta !} \p^\Ga U_j (\bz_0 ) \p_z^\beta \BN (\bx, \bz_0 ) M_{\Ga\beta}^j(B)
 + \BE(\bx),
\eeq
where the error term satisfies
\beq\label{errorbound}
\| \BE \|_{L^\infty(\p\GO)} \le C \Ge^{2d}
\eeq
for some constant $C$ independent of $\Gm$ and $\Gl$ as long as $\Gl \le \GL$ for some constant $ \GL $.
\end{thm}

We emphasize that \eqnref{asympunif} contains not only the
leading order ($\Ge^d$) term but also higher order terms up to
$\Ge^{2d-1}$. The terms higher than $\Ge^{2d}$ are expressed in
terms of not only EMTs but also interactions between the boundary
$\p\GO$ and the inclusion(s), and become much more complicated.

\medskip

\noindent{\sl Proof of Theorem \ref{firstthm}}.
We closely follow the proof in \cite{AKNT02}. By \eqnref{reprek-bvp}, the solution $\Bu$ can be written in the form:
\beq  \label{uk_repres}
\Bu(\bx) =\begin{cases}
 \Bh(\bx) +\mathcal{S}_D^0[\vp](\bx),  \quad &\bx \in \GO\setminus \overline{D}, \\
\Scal_D [\psi](\bx), \quad &\bx \in D,
\end{cases}
\eeq
where $\Bh$ is the function given by \eqnref{defBh} and
$\vp,\psi \in H^{-1/2}(\p D)$ are the solutions to
\beq \label{int_eq2}
\begin{cases}
\ds \Scal_D[\psi](\bx) - \mathcal{S}_D^0[\vp](\bx) = \Bh(\bx) \\
\nm
\ds \pd{\Scal_D[\psi]}{\Gn}\Big|_{-}(\bx) - \pd{\mathcal{S}_D^0 [\vp]}{\Gn_0}\Big|_{+}(\bx) = \pd{\Bh}{\Gn_0}(\bx)
\end{cases}
\mbox{for } \bx \in \p D.
\eeq

Let $\tilde{\varphi}(\bx): =\Ge \vp (\Bz_0 + \Ge \bx)$ and
$\tilde{\psi}(\bx): =\Ge \psi (\Bz_0 + \Ge \bx)$. By a change of
variables, \eqnref{int_eq2} can be scaled into, for $\bx \in \p
B$, \beq \label{int_eq3}
\begin{cases}
\ds \Scal_B[\tilde{\psi}](\bx) - \Scal_B^0[\tilde{\vp}](\bx) = \Bh(\Bz_0 + \Ge \bx)
-\delta_{d2} \frac{\Ga\ln\Ge}{2\pi}\int_{\p B} \tilde{\psi} , \\
\nm
\ds \pd{\Scal_B[\tilde{\psi}]}{\Gn}\Big|_{-}(\bx) - \pd{\Scal_B^0[\tilde{\vp}]}{\Gn_0}\Big|_{+}(\bx)= \pd{}{\Gn_0}\Bh(\Bz_0 + \Ge \bx),
\end{cases}
\eeq where $\delta_{d2}$ is the Kronecker delta function. By the
Taylor expansion of $\Bh=(h_1, \ldots, h_d)$, we have \beq
\label{hexp} \Bh(\Bz_0 + \Ge \bx)=  \sum_{j=1}^d\sum_{|\Ga|=0}^d
\frac{\Ge^{|\Ga|}}{\Ga !} \p^\Ga h_j(\bz_0) \bx^\Ga\Be_j   +
O(\Ge^{d+1}). \eeq Here the error term is independent of
$\Gl,~\mu$, and $\p^\Ga\Bh(\bz_0)$ is bounded independently of
$\Gl,~\mu$ as long as $\Gl$ is bounded because of \eqnref{hbound}.
By \eqnref{hexp} and the linearity of \eqnref{int_eq3} we have
\beq \label{vpexp} \tilde{\vp} (\bx)= \sum_{j=1}^d
\sum_{|\Ga|=1}^d \frac{\Ge^{|\Ga|}}{\Ga !} \p^\Ga h_j(\bz_0)
\vp_\Ga^j   + O(\Ge^{d+1}), \eeq where  $\vp_\Ga^j$ is given in
the definition of EMT. Here we see from Lemma \ref{potential} that
the error term in \eqnref{vpexp} is uniform with respect to
$\lambda,~\mu$ as long as $\Gl$ is bounded.

It is known that
\beq
\bu(\bx) = \BU(\bx) - \int_{\p D} \BN(\bx,\by) \vp(\by) d\sigma(\by), \quad \bx \in \p\GO
\eeq
(see \cite{AKNT02}).
Using \eqnref{vpexp} and the Taylor expansion of $\BN(\bx,\by)$, we have
$$ \BN(\bx, \bz_0 +\Ge \by) = \sum_{|\beta|=0}^\infty \frac{1}{\beta!} \Ge^{|\beta|} \p_\bz^\beta \BN(\bx,\bz_0) \by^\beta,\quad \bx\in \p \GO, $$
and so, for $\bx\in \p \GO$,
\begin{align}
\bu(\bx)& = \BU(\bx) - \int_{\p B} \BN(\bx, \bz_0 + \Ge\by)\tilde{ \vp}(\by)\Ge^{d-2} d\sigma(\by)\nonumber\\
&=  \BU(\bx)-  \sum_{j=1}^d \sum_{|\Ga|=1}^d  \sum_{|\beta|=1}^{d+1-|\Ga|} \frac{\Ge^{|\Ga|+|\beta|+d-2}}{\Ga !\beta !} \p^\Ga h_j (\bz_0 ) \p_z^\beta \BN (\bx, \bz_0 ) M_{\Ga\beta}^j
 + O(\Ge^{2d}).  \label{UHexp}
\end{align}

The formula \eqnref{UHexp} implies in particular that
$$ \| \Bu -  \BU \|_{H^{1/2}(\p \GO)} = O(\Ge^d) ,$$
where $O(\Ge^d)$ is uniform with respect to $\Gl$ and $\mu$.
Since
\begin{align*}
 \Bh(\bx)&= \int_{\p\GO} \pd{{\bf \Gamma^0} (\bx-\by)}{\nu_0(\by)} \Bu (\by) \, d \sigma (\by) - \int_{\p\GO} {\bf \Gamma^0} (\bx-\by) \Bg(\by) \, d \sigma (\by)\\
&= \BU(\bx) +\int_{\p\GO} \pd{{\bf \Gamma^0} (\bx-\by)}{\nu_0(\by)}( \Bu -\BU) (\by) \, d \sigma (\by),
\end{align*}
we have
\beq
 \p^\Ga \Bh (\bz_0 )=  \p^\Ga \BU (\bz_0 ) + O(\Ge^d),
\eeq
independently of $\Gl$ and $\mu$.  By substituting this into \eqnref{UHexp},  we have \eqnref{firstthm}. \qed

\section{The case of multiple inclusions}

So far we deal with the case where the inclusion $D$ is a simply connected inclusion, but the methods and results of this paper work even when there are multiple simply connected inclusions.

Let us make a brief remark on the case when $D$ has $n$ disjoint simply connected components, say $D_1, \ldots, D_n$. In this case, The solution $\Bu$ to \eqnref{trans} is represented as
\beq\label{reprek2}
\Bu(\bx) =\begin{cases}
\ds \Bh(\bx) + \sum_{j=1}^n \Scal^0_j [\vp^{(j)}](\bx),  \quad &\bx \in D^c , \\
\nm
\Scal_j [\psi^{(j)}](\bx), \quad &\bx \in D_j, \ j=1, \ldots, n,
\end{cases}
\eeq
where $\Scal^0_j$ and $\Scal_j$ denote single layer potentials on $\p D_j$, and $\vp^{(j)},\psi^{(j)}$ are the solutions to
\beq
\begin{cases}
\ds \Scal_j[\psi^{(j)}](\bx) - \sum_{j=1}^n \Scal^0_j [\vp^{(j)}](\bx) = \Bh(\bx) \\
\nm
\ds \pd{\Scal_j[\psi^{(j)}]}{\Gn}\Big|_{-}(\bx) - \sum_{j=1}^n \pd{\Scal^0_j [\vp^{(j)}]}{\Gn_0}\Big|_{+}(\bx) = \pd{\Bh}{\Gn_0}(\bx)
\end{cases}
\mbox{for } \bx \in \p D_j, \ j=1, \ldots, n.
\eeq
If $\Gk$ is either $0$ or $\infty$, we have a similar representation:
\beq\label{repreinfty2}
\Bu_\Gk(\bx) = \Bh(\bx) + \sum_{j=1}^n \Scal^0_j [\vp_{\Gk}^{(j)}](\bx),  \quad \bx \in D^c ,  \quad \Gk=0, \infty,
\eeq
where $\vp_\Gk^{(j)}$ satisfies appropriate integral equations. We can show in a similar way that
\beq
\sum_{j=1}^n \| \vp^{(j)} - \vp_\infty^{(j)} \|_{H^{-1/2}(\p D_j)} \le \frac{C}{\sqrt{\Gm}} \sum_{j=1}^n \left \| \pd{\Bh}{\Gn_0} \right\|_{H^{-1/2}(\p D_j)}
\eeq
and
\beq
\sum_{j=1}^n \| \vp^{(j)} - \vp_0^{(j)} \|_{H^{-1/2}(\p D_j)} \le C(\Gk+\mu)^{1/4} \sum_{j=1}^n \left \| \pd{\Bh}{\Gn_0} \right\|_{H^{-1/2}(\p D_j)}.
\eeq
So, Theorem \ref{thm1} and Theorem \ref{firstthm} are valid even if $D$ has several components.

It is worth mentioning that if there are multiple inclusions, the convergence depends on the distance between inclusions since $|\nabla \Bu|$ may be arbitrarily large as the distance between inclusions tends to zero. This fact was proved in, for example, \cite{AKL, Y, BLY1} for the conductivity problem. (See \cite{ACKLY} for an extensive list of recent papers on this problem.) For the elasticity problem, it is shown in \cite{KK} by numerical computations that $|\nabla \Bu|$ may blow up as the distance between inclusions tends to zero.

\section{Conclusion}

In this paper we have used layer potential techniques to prove
uniform convergence in an appropriate function spaces of solutions
to the Lam\'e system as the bulk and shear moduli tend to extreme
values (zero or $\infty$) provided that the compressional modulus
is bounded. Making use of this result, we have shown that the
asymptotic expansion of the solution due to the presence of
diametrically small inclusions is uniform with respect to the bulk
and shear moduli. These results are obtained under the assumption
that the Lam\'e parameters of the background and inclusions are
constant. We expect that the same results hold even if the
background Lam\'e parameters are not constants, but variables,
even though the methods of this paper do not apply to that case.
Another interesting case is when the inclusion is thin, and the
thickness tends to zero \cite{BF}. In this case we expect that the
asymptotic expansion is not uniform.

\end{document}